\documentclass [12pt,a4paper,reqno]{amsart}
\input amssymb.sty
\usepackage{amsmath,amsthm}
\usepackage{amsfonts}
\usepackage{amssymb}
\allowdisplaybreaks

\newcommand{\be}{\begin{equation}}
\newcommand{\ef}{\end{equation}}
\chardef\bslash=`\\ % p. 424, TeXbook
\newtheorem{thm}{Theorem}[section]
\newtheorem*{thm*}{Theorem}

\newtheorem{lem}[thm]{Lemma}
\newtheorem{prop}[thm]{Proposition}

\theoremstyle{definition}

\newtheorem*{remark*}{Remarks}
\newtheorem*{defn*}{Definition}

\theoremstyle{remark}

\numberwithin{equation}{section}

\newcommand{\G}{\Gamma}
\newcommand{\wt}{\widetilde}
\newcommand{\wh}{\widehat}
\newcommand{\fc}{\frac}
\newcommand{\bk}{\bigskip}
\newcommand{\iy}{\infty}
 \renewcommand{\sectionmark}[1]{}

\newcommand{\Be}{Beltrami}
\newcommand{\Gr} {Grunsky}
\newcommand{\hol} {holomorphic}
\newcommand{\qc} {quasiconformal}
\newcommand{\sh} {subharmonic}
\newcommand{\psh} {plurisubharmonic}

\newcommand{\ve}{\varepsilon}

\newcommand{\Te} {Teichm\"{u}ller}
\newcommand{\Ko} {Kobayashi}
\newcommand{\Ca} {Carath\'{e}odory}
\newcommand{\uTs} {universal Teichm\"{u}ller space}

 \usepackage{amsfonts}
\newcommand{\field}[1]{\mathbb{#1}}

\newcommand{\D}{\Delta}
\newcommand{\om}{\omega}
\newcommand{\z}{\zeta}
\newcommand{\ov}{\overline}
\newcommand{\vp}{\varphi}
\newcommand{\hC}{\wh{\field{C}}}
\newcommand{\C}{\field{C}}

\newcommand{\B}{\mathbf{B}}
\newcommand{\T}{\mathbf{T}}
\newcommand{\Belt}{\mathbf{Belt}}
\newcommand{\Hol}{\operatorname{Hol}}

\newcommand{\vk} {\varkappa}
\newcommand{\x} {\mathbf x}
\renewcommand{\a} {\alpha}

\newcommand{\ld}{\lambda}
\newcommand{\kp}{\kappa}

\begin{document}

\title{Hyperbolic metrics, homogeneous holomorphic
functionals and Zalcman's conjecture}
\author{Samuel L. Krushkal}

\begin{abstract}
We show, using the \Ko \ and \Ca \ metrics on special  
holomorphic disks in the \uTs,
that a wide class of \hol \ functionals on
the space of univalent functions in the disk is maximized by
the Koebe function or by its root transforms; their extremality
is forced by hyperbolic features.
As consequences, this implies the proofs of
the famous Zalcman and Bieberbach conjectures.
\end{abstract}

\date{\today\hskip4mm({hypZB2.tex})}

\maketitle

\bigskip

{\small {\textbf {2010 Mathematics Subject Classification:}
Primary: 30C50, 30F60, 32Q45; Secondary 30C55, 30C62}

\medskip

}

\bigskip

\markboth{Samuel L. Krushkal}{Hyperbolic metrics,
homogenous functionals and Zalcman's conjecture}
\pagestyle{headings}

\section{Introduction. Main theorems}

In this paper, we build a bridge linking the hyperbolic
\Ko \ and \Ca \ metrics of the \uTs \ and \Gr \ inequalities
with two famous conjectures in classical geometric complex
analysis and prove these conjectures.

\subsection{General homogeneous \hol \ functionals} 
The  \hol \ functionals on
the classes of univalent functions depending on the Taylor
coefficients of these functions play an important role
in various geometric and physical applications of
complex analysis, for example, in view of their connection with
string theory and with a \hol \ extension of the Virasoro algebra.
These coefficients reflect the fundamental intrinsic features
of conformal maps. Thus estimating them still remains an
important problem in geometric function theory.

We consider the univalent functions on the unit disk
$\D = \{|z| < 1\}$ normalized by
$$
f(z) = z + \sum\limits_{n=2}^\iy a_n z^n.
$$
These functions form the well-known class $S$.
Their inversions
$F_f(z) = 1/f(1/z)$ belong to the class $\Sigma$ of univalent
nonvanishing functions
\be\label{1.1}
F(z) = z + b_0 + b_1 z^{-1} + b_2 z^{-2} + \dots \
\quad (F(z) \ne 0)
\end{equation}
on the complementary disk
$\D^* = \{z \in \hC = \C \cup \{\iy\}: \ |z| > 1\}$.
Easy computations yield that the coefficients $a_n$ and $b_j$
are related by
\be\label{1.2}
b_0 + a_2 = 0, \quad
b_n + \sum \limits_{j=1}^{n} b_{n-j} a_{j+1} + a_{n+2} = 0,
\quad n = 1, 2, ... \ ,
\end{equation}
which implies successively the representations
of $a_n$ by $b_j$. One gets
\be\label{1.3}
a_n = (- 1)^{n-1} b_0^{n-1} - (- 1)^{n-1}  (n - 2) b_1 b_0^{n-3} +
\text{lower terms with respect to} \ b_0;
\end{equation}
in particular,
$$
\begin{aligned}
a_2 &= - b_0, \ \  a_3 = - b_1 + b_0^2, \ \
a_4 = - b_2 + 2 b_1 b_0 - b_0^3,  \\
a_5 &= - b_3 + 2 b_2 b_0 + b_1^2 - 3 b_1 b_0^2 + b_0^4, \\
a_6 &= - b_4 + 2 b_3 b_0 + 2 b_2 b_1 - 3 b_2 b_0^2 - 3 b_1^2 b_0
+ 4 b_1 b_0^3 - b_0^5, \\
a_7 &= b_0^6 - 5 b_1 b_0^4  - b_1 ^3 + 4 b_2 b_0^3 + b_2^2
+ (6 b_1^2 - 3 b_3) b_0^2 \\
&+ 2 b_1 b_3 + (-6 b_1 b_2 + 2 b_4) b_0 - b_5, \ \dots
\end{aligned}
$$
We shall essentially use this connection.

Consider a general \hol \  distortion functional on $S$ of the form
\be\label{1.4}
J(f) = J(a_2, \dots \ , a_n;
(f^{(\a_1)}(z_1)); \dots \ ; (f^{(\a_p)}(z_p))),
\end{equation}
where $z_1, \dots \ , z_p$ are the distinct fixed points in
$\D \setminus \{0\}$ with assigned orders
$m_1, \ \dots \ , m_p$, respectively,
$(f^{(\a_1)}(z_1)) = f^{\prime\prime}(z_1), \dots, f^{(m_1)}(z_1); \
(f^{(\a_p)}(z_p)) = f^{\prime\prime}(z_p), \dots, f^{(m_p)}(z_p)$.
Assume that $J$ is a polynomial in all of its variables.

Substituting the expressions of $a_j$ by $b_m$ from
(1.2) and calculating $f^{(q)}(z_j)$ in terms of $F_f$,
one obtains a polynomial $\wt J(F)$ of the Taylor coefficients
$b_0, b_1, \dots \ , b_{n-2}$ and of the corresponding
derivatives $F_f^{(q)}(\z_j)$ at the points
$\z_j = 1/z_j \in \D^* \setminus \{\iy\}$,
regarded as a representation of $J(f)$ on the
class $\Sigma$.
Here $q = 2, \ \dots \ , m_j, \ j = 1, \ \dots \ , p$.

Assume that the functional (1.4) is homogeneous with a degree
$d = d(J)$ (depending on $n$ and $m_1, \ , \dots, \ , m_p$)
with respect to homotopy
$$
f(z, t) = t^{-1} f(t z) = z + a_2 t + a_3 t^2 + \dots: \
\D \times \ov \D \to \C
$$
such that $f(z, 0) \equiv z, \ f(z, 1) = f(z)$ so that
$$
J(f_t) = t^d J(f).
$$
This homotopy is a special case of \hol \ motions with
complex parameter $t$ running over the disk $\D$.
The functional $\wt J(F)$ on $\Sigma$ admits a
similar homogeneity.

The existence of extremal functions of $J(f)$ and $\wt J(F)$
follows from compactness of both classes $S$ and $\Sigma$ in
the topology of locally uniform convergence on $\D$ and $\D^*$,
respectively.

\subsection{The Zalcman conjecture}
There were several classical conjectures about these coefficients. 
They include the Bieberbach conjecture that in the class $S$  the
coefficients are estimated by $|a_n| \leq n$, as well as several other 
well-known conjectures that imply the Bieberbach conjecture.
Most of them have been proved by the de Branges theorem \cite {DB}.

In the 1960s, Lawrence Zalcman posed the conjecture that for
any $f \in S$ and all $n \ge 3$,
\be\label{1.5}
|a_n^2 - a_{2n-1}| \le (n-1)^2,
\end{equation}
with equality only for the {\bf Koebe function}
\be\label{1.6}
\kp_\theta(z) =
\fc{z}{(1 - e^{i \theta} z)^2} = z + \sum\limits_2^\iy
n e^{- i(n-1) \theta} z^n, \quad 0 \le \theta \le 2 \pi,
\end{equation}
which maps the unit disk onto the complement of the ray
$$
w = -t e^{-i \theta}, \ \ \fc{1}{4} \le t \le \iy.
$$
This remarkable conjecture also implies the Bieberbach conjecture.
This still is an intriguing very difficult
open problem for all $n > 6$.

The original aim of Zalcman's conjecture was
to prove the Bieberbach conjecture using the famous Hayman theorem
on the asymptotic growth of coefficients of individual functions,
which states that {\em for each
$f \in S$, we have the inequality
$$
\lim\limits_{n \to \iy} \fc{|a_n|}{n} = \a \le 1,
$$
with equality only when $f = \kp_\theta$; here}
$\a = \lim\limits_{r \to 1} (1 - r)^2 \max_{|z|=r} |f(z)|$
(see \cite{Ha}).

Indeed, assuming that $n$ is sufficiently large and estimating
$a_{2n-1}$ in (1.5) by $|a_{2n-1}| \le 2n -1$, one obtains
$$
|a_n|^2 \leq (n -1)^2 + |a_{2n -1}| \leq (n - 1)^2 + 2n - 1 = n^2,
$$
which proves the Bieberbach conjecture for this $n$, and
successively for all preceding coefficients.

It was realized almost immediately that the Zalcman conjecture implies
the Bieberbach conjecture, and in a very simple fashion, without
Hayman's result and without other prior results from the theory of
univalent functions.

Note that the case $n = 2$ is rather simple and somewhat exceptional.
The inequality
$|a_2^2 - a_3| \le 1$
is well known, but in this case there are two
extremal functions of different kinds: the Koebe function $\kp_\theta(z)$
and the odd function
\be\label{1.7}
\kp_{2,\theta}(z) := \sqrt{\kp_\theta(z^2)}
= \sum\limits_{n=0}^\iy e^{i n \theta} z^{2n+1}.
\end{equation}
The estimate (1.5) was established for $n = 3$ in \cite{Kr4}
and recently for $n = 4, 5, 6$ in \cite{Kr8}
(without uniqueness of the extremal function).
In \cite{BT}, \cite{Ma}, this conjecture was proved for
certain special subclasses of $S$.

\subsection{Main theorems}
It is well known that the Koebe function $\kappa_\theta$
is extremal for many variational problems in the theory
of conformal maps (accordingly, its root transforms
\be\label{1.8}
\kp_{m,\theta}(z) = \kp_\theta(z^m)^{1/m}
= \fc{z}{(1 - e^{i \theta} z^m)^{2/m}}
= z + \fc{2 e^{i \theta}}{m} z^{m+1} + \fc{m - 2}{m^2} z^{2m+1}
+ \dots, \quad m = 2, 3, \dots,
\end{equation}
are extremal among the maps with symmetries).

Our first main theorem sheds new light on this phenomenon and
provides a large class of functionals maximized by these
functions.

\begin{thm}
Let $J(f)$ be a homogeneous polynomial functional
on $S$ of the form (1.4) whose representation $\wt J(F_f)$
in the class $\Sigma$ does not contain free terms
$c_d b_0^d$ but contains nonzero terms with the coefficient
$b_1$ of inversions $F_f$.
Then for all $f \in S$, we have the sharp bound
\be\label{1.9}
|J(f)| \le \max_m |J(\kappa_{m,\theta})|,
\end{equation}
and this maximum is attained on some
$\kappa_{m_0,\theta} \ (m_0 \ge 1)$.
If $J$ has an extremal with
\be\label{1.10}
b_1 = a_2^2 - a_3 \ne 0,
\end{equation}
then $|b_1| = 1$ and
\be\label{1.11}
|J(f)| \le \max \{|J(\kappa_\theta)|, |J(\kappa_{2,\theta})|\}.
\end{equation}
\end{thm}

The assumption (1.10) is equivalent to
$$
S_f(0) = - \lim\limits_{z\to \iy} z^4 S_{F_f}(z) \ne 0.
$$
The examples of some well-known functionals, for example,
$J(f) = a_2^2 - \a a_3$ with $0 < \a < 1$ and
$J(F_f) = b_m \ (m > 1)$, show that the assumptions on
the initial coefficients $b_0$ and $b_1$  cannot be omitted.

\bk
The Zalcman functional
$$
J_n(f) = a_n^2 - a_{2n-1}
$$
is a special case of (1.4) with homogeneity degree $2n - 2$.
For this functional, we obtain from Theorem 1.1 a
complete result proving the Zalcman conjecture.

\begin{thm} For all $f \in S$ and any $n \ge 3$,
we have the sharp estimate (1.5), with equality only for
$f = \kappa_\theta$.
\end{thm}

As a consequence, one obtains also a new proof of
the Bieberbach conjecture.

\bk
Theorem 1.1 also provides other new distortion theorems
concerning the higher coefficients. These results are
presented in the last section.

\subsection{}
It suffices to find the bound of $J$ on functions admitting
\qc \ extensions across the unit circle and make the closure of
this set in weak topology determined by locally uniform
convergence on $\D$. Such functions are naturally connected with
the \uTs \ $\T$. The original functional $J(f)$ is lifted
to a \hol \ functional on a fiber space over $\T$, but its
growth is controlled by hyperbolic metrics on the base space $\T$.
Extremity of the Koebe function or of its root
transforms is intrinsically connected with the features
of these metrics on appropriate geodesic disks.

\section{Background}

We briefly present here certain results underlying the proof of
Theorem 1.1. This exposition is adapted to our special cases.

\subsection{Universal \Te \ space}
We shall use some deep geometric and potential results concerning
the \Te \ space of the disk, called also the \uTs.

{\em (a)} \ The \uTs \ $\T$ is the space of quasisymmetric
homeomorphisms
of the unit circle $S^1 = \partial \D$ factorized by M\"{o}bius maps.
The canonical complex Banach structure on $\T$ is defined by
factorization of the ball of the {\bf \Be \ coefficients} (or complex
dilatations)
\be\label{2.1}
\Belt(\D)_1 = \{\mu \in L_\iy(\C): \ \mu|\D^* = 0, \ \|\mu\| < 1\},
\end{equation}
letting $\mu_1, \mu_2 \in \Belt(\D)_1$ be equivalent if the
corresponding \qc \ maps $w^{\mu_1}, w^{\mu_2}$
(solutions to the \Be \ equation
$\partial_{\ov{z}} w = \mu \partial_z w$ with $\mu = \mu_1, \mu_2$)
coincide on the unit circle $S^1 = \partial \D^*$ (hence, on $\ov{\D^*}$).
The equivalence classes $[w^\mu]$ are in one-to-one correspondence with
the Schwarzian derivatives
$$
S_w(z) := \Bigl(\fc{w^{\prime\prime}}{w^\prime}\Bigr)^\prime
- \fc{1}{2} \Bigl(\fc{w^{\prime\prime}}{w^\prime}\Bigr)^2
= \fc{w^{\prime\prime\prime}}{w^\prime} -
\fc{3}{2} \Bigl(\fc{w^{\prime\prime}}{w^\prime}\Bigr)^2, \quad
w = w^\mu|\D^*.
$$

Note that for each locally univalent function $w(z)$ on a simply connected
hyperbolic domain $D \subset \hC$ its Schwarzian derivative $S_w$
belongs to the complex Banach space $\B(D)$ of
hyperbolically bounded \hol \ functions on $D$ with the norm
$$
\|\vp\|_\B = \sup_D \ld_D^{-2}(z) |\vp(z)|,
$$
where $\ld_D(z) |dz|$ is the hyperbolic metric on $D$ of Gaussian
curvature $- 4$; hence $\vp(z) = O(z^{-4})$ as $z \to \iy$
if $\iy \in D$. In particular, for $D = \D$,
\be\label{2.2}
\ld_\D(z) = 1/(1 - |z|^2).
\end{equation}
The space $\B(D)$ is dual to the Bergman space $A_1(D)$,
a subspace of $L_1(D)$ formed by integrable \hol \ functions
on $D$.

The derivatives $S_{w^\mu}(z)$ with $\mu \in \Belt(\D)_1$
range over a bounded domain in the space $\B = \B(\D^*)$.
This domain models the \uTs \ $\T$, and the defining projection
$$
\phi_\T(\mu) = S_{w^\mu}: \ \Belt(\D)_1 \to \T
$$
is a \hol \ map from $L_\iy(\D)$ to $\B$. This map is a split submersion,
which means that $\phi_\T$ has local \hol \ sections (see, e.g., [GL]).

The above definition of $\T$ requires a complete normalization of
maps $w^\mu$, which uniquely define the values of $w^\mu$ on $\D^*$
by their Schwarzians.  We shall use the condition $w^\mu(0) = 0$
going over from $w^\mu$ to the maps
$$
w_1^\mu(z) = w^\mu(z) - w^\mu(0)
= z - \fc{1}{\pi} \iint\limits_\D
\fc{\partial w^\mu}{\partial \ov \z}
\left( \fc{1}{\z - z} - \fc{1}{\z} \right) d \xi d \eta
\quad (\z = \xi + i \eta),
$$
which does not reflect on the Schwarzians. We identify the
bounded domain in $\B$ filled by $S_{w_1^\mu}$ with the
space $\T$.

\bk
The intrinsic {\bf \Te \ metric} of the space $\T$ is defined by
\be\label{2.3}
\tau_\T (\phi_\T (\mu), \phi_\T
(\nu)) = \frac{1}{2} \inf \bigl\{ \log K \bigl( w^{\mu_*} \circ
\bigl(w^{\nu_*} \bigr)^{-1} \bigr) : \ \mu_* \in \phi_\T(\mu),
\nu_* \in \phi_\T(\nu) \bigr\};
\end{equation}
it is generated by the Finsler structure
\be\label{2.4}
F_\T(\phi_\T(\mu), \phi_\T^\prime(\mu) \nu) =
\inf \{\|\nu_*/(1 - |\mu|^2)^{-1}\|_\iy:
\phi_\T^\prime(\mu) \nu_* = \phi_\T^\prime(\mu) \nu\}
\end{equation}
on the tangent bundle $\mathcal T(\T) = \T \times \B$ of $\T$
(here $\mu \in \Belt(\D)_1$ and $\nu, \nu_* \in  L_\iy(\C)$).
This structure is locally Lipschitz (cf. \cite{EE}).

The smallest dilatation $k(F) = \inf \|\mu_{\wh F}\|_\iy$ among
\qc \ extensions of $F|\D^*$ onto $\hC$ is called the
\textbf{\Te \ norm} (or dilatation) of $F$.

The space $\T$ as a complex Banach manifold also has invariant metrics
(with respect to its biholomorphic automorphisms);
the largest and the smallest invariant metrics are
the \Ko \ and the \Ca \ metrics, respectively.
Namely, the \textbf{\Ko \ metric} $d_\T$ on $\T$ is the
largest pseudometric $d$ on $\T$ which does not get increased by \hol \ maps
$h: \ \D \to \T$ so that for any two points $\psi_1, \ \psi_2 \in \T$,
$$
d_\T(\psi_1, \psi_2) \leq \inf \{d_\D(0,t):
\ h(0) = \psi_1, \ h(t) = \psi_2\},
$$
where $d_\D$ is the hyperbolic metric on $\D$ with the
differential form (2.2).
This distance is connected with the \Te \ norm of $f$ by
$k(f) = \tanh d_\T(\mathbf 0, S_f)$.

The {\bf \Ca \ distance} between $\psi_1$ and $\psi_2$ in $\T$ is
$$
c_\T(\psi_1, \psi_2) = \sup d_\D(h(\psi_1), h(\psi_2)),
$$
where the supremum is taken over all  \hol \ maps $h: \ \T \to \D$.

The \Ko \ metric is the integrated form of its infinitesimal  
Finsler metric defined for the points $(\psi, v) \in \mathcal T(\T)$ by
$$
\mathcal K_\T(\psi, v) = \inf \{1/r: \ r > 0, \ h \in \Hol(\D_r, \T),
h(0) = \psi, h^\prime (0) = v\}, 
$$
where $\Hol(\D_r, \T)$ denotes the collection of \hol \ maps of 
the disk $\D_r = \{|z| < r\}$ into $\T$. 
For the general properties of invariant metrics
we refer, for example, to \cite{Di}, \cite{Ko}.

The Royden-Gardiner theorem states that the \Ko \ and
\Te \ metrics are equal on all \Te \ spaces (cf. [EKK], [EM],
[GL], [Ro]). This fundamental fact underlies many applications of
the \Te \ space theory.

\bk
{\em (2)} \ A strengthened version of the Royden-Gardiner theorem
for the \uTs \ is given by

\begin{prop} \cite{Kr4}. The infinitesimal \Ko \ metric
$\mathcal K_\T(\psi, v)$ on the tangent bundle $\mathcal T(\T)$ of the
\uTs \ $\T$ is logarithmically \psh \ in $\psi \in \T$, equals the
canonical Finsler structure $F_\T(\psi, v)$ on $\mathcal T(\T)$
generating the \Te \ metric of $\T$ and has constant \hol \ sectional
curvature $- 4$.
\end{prop}

It implies that {\em the \Te \ distance $\tau_\T(\vp, \psi)$
is logarithmically \psh \ in each of its variables and hence
the pluricomplex Green function of the space $\T$ equals
\be\label{2.5}
g_\T(\vp, \psi) = \log \tanh \tau_\T(\vp, \psi) = \log k(\vp, \psi),
\end{equation}
where $k(\vp, \psi)$ denotes the extremal dilatation of \qc \ maps
determining
the \Te \ distance between the points $\vp$ and $\psi$ in $\T$. }

Recall that the {\bf pluricomplex Green function} $g_D(x, y)$
of a domain $D$ in a complex Banach manifold $X$ with pole $y$ is
defined by
$g_D(x, y) = \sup u_y(x) \ (x, y \in D)$ 
and followed by upper regularization
$v^*(x) = \limsup_{x^\prime \to x} v(x^\prime)$, taking the supremum 
over all \psh \ functions $u_y(x): \ D \to [-\iy, 0)$ such that
$$
u_y(x) = \log \|x - y\|_X + O(1)
$$
in a neighborhood of the pole $y$. Here $\|\cdot\|_X$ denotes the
norm on the space modeling $X$, and the remainder term $O(1)$
is bounded from above.
The Green function $g_D(x, y)$ is a maximal \psh \ function on
$D \setminus \{y\}$ (unless it is identically $- \iy$).

\bk
{\em (c)} \ The assumption $F^\mu(0) = 0$ for \qc \ extensions
of a function $F \in \Sigma^0$ made above
ensures nonvanishing $F$ in $\D^*$ (and $f(z) = 1/F(1/z) \in S$).

In addition, completely normalized maps $F^\mu(z)$
are \hol \ functions of their \Be \ coefficients
$\mu \in \Belt(\D)_1$ as well as of their Schwarzians.
The same holds for the coefficients $b_m$ of $F^\mu$.

The Schwarzian equations $S_w = \vp$ for $F \in \Sigma$ and
for its inversion $f \in S$ and the \Be \ equation for
\qc \ extensions of these functions determine their solutions
up to linear transformations, which for $F$ are
the translations $w \mapsto w + b_0$ and for $f$
have the form
$$
w \mapsto w/(1 - \a w) = w + \a w^2 + \cdots
$$
with $\a$ determined by the coefficients $a_2$.
The admissible values of $b_0 = - a_2$, which are the free
terms of corresponding  $F_1 \in \Sigma$ having the same
Schwarzian, are only those which range over the closed
domain $\hC \setminus F(\D^*)$.

Since the original normalization of the functions
from $\Sigma$ and $S$  includes only two conditions,
the initial functional $J_n$ must be considered
on the \textbf{fiber space} $\mathcal F(\T)$ over $\T$, which
is modeled as a bounded domain in the space
$\B \times \D(0, 2)$, whose points are the pairs $(\vp, a)$
where $\vp$ are the Schwarzians of $F \in \Sigma^0$
with $F(0) = 0$
and $a$ are equal to the second coefficients $a_2$ of their
inversions in $S$.
The defining projection of this space
$$
\pi_{\mathcal F}: \ (S_{F^{\wt \mu}}, a_2^\mu) \to S_{F^{\wt \mu}}
$$
is a \hol \ split submersion.
The fibers $\pi_{\mathcal F}^{-1}(\vp)$ over the base points $\vp = S_F$
coincide with the complementary domains $\hC \setminus \ov{F(\D^*)}$,
giving the admissible values of $a_2$.

Note also that the space $\mathcal F(\T)$ is biholomorphically
isomorphic to the Bers fiber space over $\T$ (both fiber spaces
have the same base and differ only by an additional normalization
of maps $F \in \Sigma$ with a given Schwarzian $\vp = S_F$)
and thus is isomorphic to the \Te \ space of the punctured disk
$\D \setminus \{0\}$ (cf. \cite{Be2}).

\subsection{The Grunsky operator}
The complex geometry of the \uTs \ is closely connected with
the \Gr \ inequalities technique which arose from
investigating the univalence problem in \cite{Gr}.

Any function $F \in \Sigma$ (and similarly any $f \in S$)
determines its \Gr \ operator (matrix)
$\mathcal G_F = (\a_{mn}(F))$, where the
\Gr \ coefficients $\a_{mn}$ are determined by the
expansion
$$
\log \fc{F(z) - F(\z)}{z - \z} = - \sum\limits_{m,n=1}^\iy
\a_{mn} z^{-m} \z^{-n}, \quad (z, \z) \in (\D^*)^2,
$$
with the principal branch of logarithmic function,
satisfy the inequalities
\be\label{2.6}
\big\vert \sum\limits_{m,n=1}^\iy \ \sqrt{mn} \ \a_{mn} x_m x_n
\big\vert \le  k.
\end{equation}
Here $\mathbf x = (x_n)$ runs over the unit sphere $S(l^2)$
of the Hilbert space $l^2$ with norm
$\|\x\| = \Bigl(\sum\limits_1^\iy |x_n|^2 \Bigr)^{1/2}$,
and $k = k(F) \le 1$ is the \Te \ norm of $F$
(cf. \cite{Gr}, \cite{Ku1}).
The quantity
$$
\vk(F) = \sup \Big\{\Big\vert \sum\limits_{m,n=1}^\iy \
\sqrt{mn} \ \a_{mn} x_m x_n \Big\vert:
\mathbf x = (x_n) \in S(l^2) \Big\} \le 1 
$$
is called the \textbf{\Gr \ norm} of $F$.
It equals the norm of $\mathcal G_F$ regarded as a linear
operator $l^2 \to l^2$.

The functions with $\vk(F) = k(F)$ play a crucial role
in applications of \Gr \ inequalities; however, 
the set of $S_F$, on which $\vk(F) < k(F)$, is open and dense
in $\T$.
One of the underlying facts in applications is the following result.

\begin{prop} The equality $\vk(F) = k(F)$ for $f \in \Sigma^0$
holds if and only if
the function $F$ is the restriction to $\D^*$ of a
\qc \ self-map $w^{\mu_0}$ of $\hC$ with \Be \ coefficient
$\mu_0$ satisfying the condition
$\sup |\langle \mu_0, \psi\rangle_\D| = \|\mu_0\|_\iy$, 
where the supremum is taken over \hol \ functions
$\psi \in A_1^2(\D)$ with $\|\vp\|_{A_1(\D)} = 1$,
where
$$
A_1^2 = \{\psi \in A_1(\D): \ \psi = \om^2 \ \
\text{with} \ \ \om \ \ \text{\hol \ on} \ \ \D.
$$
In addition, if the equivalence class $[F]$
contains a frame map, i.e., is a Strebel point
(see Section 2.4), then
the restriction of $\mu_0$ onto the disk $\D$
must be of the form
\be\label{2.7}
\mu_0(z) = k |\psi_0(z)|/\psi_0(z) \quad \text{with} \ \
\psi_0 \in A_1^2 .
\end{equation}
\end{prop}

\medskip
The proof of this proposition is given in \cite{Kr2}, \cite{Kr7}.
It relies on the fact that the \Gr \ coefficients $\a_{mn}(S_F)$
generate the \hol \ functions
\be\label{2.8}
h_{\mathbf x}(\vp) =
\sum\limits_{m,n=1}^\iy \ \sqrt{mn} \ \a_{mn}(\vp) x_m x_n,
\end{equation}
where $\vp = S_f$ and $\x = (x_n)$ are the points of the shere
$S(l^2)$, mapping the \uTs \ $\T$ into the unit disk $\D$.
The restrictions of these functions to the disk $\{\phi_\T(s\mu_0)\}$
determine the \Ca \ distance between the points $S_{f^{s\mu_0}}$
and the origin, which by (2.6) equals the \Te \ distance.

In a special case, when the curve $F(S^1)$ is analytic,
the equality (2.7) was obtained by a different method
in \cite{Ku2}.

In particular, {\em the maps}  
\be\label{2.9}
F_{m,t}(z) = \fc{1}{\kp_{m,t}(1/z)} =
z \left(1 - \fc{t}{z^{m+1}}\right)^{2/(m+1)}
= z - \fc{2 t}{m + 1} \fc{1}{z^m} + \dots,  \quad |t| \le 1,
\end{equation}
whose extremal extensions to $\C$ have \Be \ coefficients 
$$
\mu_{F_{m,t}}(z) = t |z|^{m-1}/z^{m-1}
\quad \text{for} \ \ |z| < 1, 
$$
{\em satisfy $\vk(F_{m,t}) = k(F_{m,t})$ for odd $m \ge 1$ 
and $\vk(F_{m,t}) < k(F_{m,t})$ for even $m \ge 0$.} 

\bk
Note that holomorphy of the functions (2.8) is a consequence of the
fact that the \Gr \ coefficients $\a_{mn}$ are polynomials of
the initial coefficients $b_1, \dots, b_{m+n-1}$ of $F$
combined with the well-known inequality (cf. [Po, p. 61]) :
for any $1 \le p \le M, \ 1 \le q \le N$,
$$
\Big\vert \sum\limits_{m=p}^M \sum\limits_{n=q}^N
\sqrt{mn} \ \a_{mn} x_m x_n \Big\vert^2 \le
\sum\limits_{m=p}^M |x_m|^2  \sum\limits_{n=q}^N |x_n|^2.
$$

We mention also that both \Te \ and Grunsky norms are continuous
logarithmically \psh \ functions on $\T$ (see \cite{Kr6},
\cite{Sh}) and that, by a theorem of Pommerenke and
Zhuravlev, any $F  \in \Sigma$ with $\vk(F) \le k < 1$ has
$k_1$-\qc \ extensions to $\hC$ with $k_1 = k_1(k) \ge k$
(see [Po]; [KK1, pp. 82-84], [Zh]).

\subsection{A holomorphic homotopy of univalent function}
Similarly to the functions in $S$, one can define
for each $F \in \Sigma$ with expansion (1.1) the complex homotopy
\be\label{2.10}
F_t (z) = t F \left( \fc{z}{t} \right) =  z + b_0 t + b_1 t^2 z^{-1} +
b_2 t^3 z^{-2} + ...: \ \D^* \times \D \to \hC
\end{equation}
so that $F_0(z) \equiv z$. This implies
$$
S_{F_t}(z) = t^{-2} S_F(t^{-1} z),
$$
and moreover, this point-wise map determines a \hol \ map
\be\label{2.11}
h_F(t) =  S_{F_t}(\cdot): \ \D \to \B
\end{equation}
(see, e.g. \cite{Kr3}).
The corresponding \textbf{homotopy disk}
\be\label{2.12}
\D(S_F) = h_F(\D) \subset \T
\end{equation}
is \hol \ at noncritical points of maps (2.11).
These disks foliate the space $\T$ (and the
set $\Sigma$).

The dilatations of the homotopy maps are estimated by

\begin{prop} \cite{Kr3} (a) \ Each homotopy map $F_t$
of $F \in \Sigma$ admits $k$-\qc \
extension to the complex sphere $\hC$ with $k \le |t|^2$.
The bound $k(F_t) \le |t|^2$ is sharp and occurs only for the maps
$$
F_{b_0,b_1} (z) = z + b_0 + b_1 z^{-1}, \quad |b_1| = 1,
$$
whose homotopy maps
\be\label{2.13}
F_{b_0,b_1t^2}(z) = z + b_0 t + b_1 t^2 z^{-1}
\end{equation}
have the affine extensions
$\wh F_{b_0,b_1t^2}(z) = z + b_0 t + b_1 t^2 \ov z$
onto $\D$.

(b) \ If $F(z) = z + b_0 + b_m z^{-m} + b_{m+1} z^{-(m+1)} + \dots \
(b_m \ne 0)$
for some integer $m > 1$, then the minimal dilatation of extensions
of $F_t$ is estimated by $k(F_t) \le |t|^{m+1}$; this bound also is sharp.
\end{prop}

In the second case,
$$
h_F(0) = h_F^\prime(0) = \dots = h_F^{(m)}(0) = \mathbf 0, \
h_F^{(m+1)}(0) \ne \mathbf 0,
$$
and due to \cite{KK2},
\be\label{2.14}
k(F_t) = \fc{m + 1}{2} |b_m| |t|^{m+1} + O(t^{m+2}), \quad t \to 0.
\end{equation}
This bound is sharp.

The simplest \hol \ disks in $\Sigma^0$ are the images of \Te \
\textbf{extremal disks}
$$
\D(\psi) = \{\phi_\T(t \mu_0): \ t \in \D\} \subset \T
$$
formed by functions $F$ whose extremal extensions onto $\D$
(with minimal dilatation) have \Be \ coefficients $t \mu_0$ 
with 
\be\label{2.15}
\mu_0(z) = |\psi(z)|/\psi(z),
\end{equation}
where $\psi$ is a \hol \ function from $L_1(\D)$.
Such an extension is unique (up to a constant factor of $\psi$
and normalization of $F$).
The \Te \ disks foliate dense subsets in $\T$ and $\Sigma^0$.
Note also that the homotopy function $F_t$ has such extremal
extension for each $t \in \D$ (cf. \cite{GL}, \cite{St}) 
and that for any function (2.9) its homotopy disk $\D(S_{F_m})$ 
is of \Te \ type.

\section{Proof of Theorem 1.1}

$(a)$ \ Assume that there exists an extremal of $J(f)$ satisfying
the assumption (1.10), and consider first the functions
$f \in S$, for which this inequality is fulfilled.
The set of the corresponding Schwarzians $S_{F_f}$
is dense in $\T$, and their maps (2.12) satisfy
$$
h_F(0) = h_F^\prime(0) = \mathbf 0, \ \
h_F^{\prime\prime}(0) \ne \mathbf 0.
$$
Using the relations (1.2), we represent $J$ as a
polynomial functional on $\Sigma$, which takes the form
\be\label{3.1}
J(f) = \wt J(F_f) = \wt J_n(b_0, b_1, \dots, b_{2n-3};
F_f^{\prime\prime}(\z_1), \dots, F_f^{(m_1)}(\z_1); \dots,
F_f^{\prime\prime}(\z_p), \dots, F_f^{(m_p)}(\z_p))
\end{equation}
(where $b_0 = - a_2$ and $\z_j = 1/z_j$) and put
$$
\wt J^0(F_f) = \fc{\wt J(F_f)}{M(J)} \quad
\text{with} \ \ M(J) = \max_S |J(f)|
$$
to have a \hol \ map $\mathcal F(\T) \to \D$.

As was mentioned above, the admissible values of $- b_0 = a_2$
for 
$F(z) = z + b_0 + b_1 z^{-1} + \dots \in \Sigma^0$ with $F(0) = 0$
range over the closed domain
$F(\ov{\D}) = \hC \setminus F(\D^*)$.
In view of the maximum principle, it suffices to use only
the boundary points of domains $F(\D^*)$.

We select on the unit circle $S^1$ a dense subset
$$
e = \{Z_1, Z_2, \dots , \ Z_m, \dots\}.
$$
Then the images
$F(Z_1), \dots, \ F(Z_m), \dots$ 
for $F \in \Sigma^0$ (with normalization indicated above)
generate a sequence of \hol \ maps
\be\label{3.2}
g_m(\vp) =
\wt J^0(- F(Z_m), b_1(\vp), \dots, \ b_{2n-3}(\vp));
\{F^{(m_j)}(\z_j(\vp))\}): \ \T \to \D, \quad m = 1, 2, \dots
\end{equation}
where $\{F^{(m_j)}(\z_j(\vp))\}$ denotes the collection
$$
F^{\prime\prime}(\z_1), \dots, F^{(m_1)}(\z_1) ; \dots; \
F^{\prime\prime}(\z_p), \dots, F^{(m_p)}(\z_p).
$$
The upper envelope
$$
\mathcal J(S_F) = \sup_m |g_m(S_f)|,
$$
followed by its upper semicontinuous regularization
$\mathcal J(\vp) = \limsup\limits_{\vp^* \to \vp} \mathcal J(\vp^*)$,
yields a logarithmically \psh \ functional $\T \to [0, 1)$
so that
$$
\sup_\T \mathcal J(S_F) = \sup_{\Sigma^0} |\wt J^0(F)|
= \max_S |J(f)|/M(J).
$$
One may assume that the degree $d$ of $J$ (hence of $\mathcal J$)
is even, replacing, if needed, this functional by its square $J^2$.

In accordance with normalization $F_f(0) = 0$ of extended maps,
we pick in (2.13)
$$
F_{0,b_1}(z) = z + b_1z^{-1}.
$$
The extremal extension of this map onto $\ov \D$ is $z + b_1 \ov z$.
Now we split every homotopy function $F_t$ of $F = F_f$ by
$$
F_t(z) = z + b_0 t + b_1 t^2 z^{-1} + b_2 t^3 z^{-2} + \dots
= F_{0,b_1t^2}(z) + h(z,t).
$$
For sufficiently small $|t|$, the remainder $h$ is estimated by
$h(z,t) = O(t^3)$ uniformly in $z$ for all $|z| \ge 1$.
Then, by the well-known properties of Schwarzians, we have
$$
S_{F_t}(z) = S_{F_{b_0,b_1t^2}}(z) + \om(z, t)
= S_{F_{0, b_1t^2}}(z) + \om(z, t),
$$
where the remainder $\om$ is uniquely determined by the chain
rule
$$
S_{w_1\circ w}(z) = (S_{w_1} \circ w) (w^\prime)^2(z) + S_w(z),
$$
and is estimated in the norm of $\B$ by
\be\label{3.3}
\|\om(\cdot,t)\|_\B = O(t^3), \quad t \to 0;
\end{equation}
this estimate is uniform for $|t| < t_0$
(cf., e.g. \cite{Be1}, \cite{Kr1}).

Then, in view of holomorphy, every map (3.2) satisfies for
small $|t|$,
\be\label{3.4}
g_m(S_{F_t}) = g_m(S_{F_{0,b_1t^2}}) + O(t^{d+1}),
\end{equation}
where the term $O(t^{d+1})$ is estimated uniformly for all $n$.
Thus
$$
\mathcal J(S_{F_t}) = \mathcal J(S_{F_{0,b_1t^2}}) + O(t^{d+1})
\quad t \to 0.
$$
Combining with the equality
$$
\mathcal J(S_{F_t}) = t^d \mathcal J(S_F) +  O(t^{d+1}),
$$
which follows from $d$-homogeneity of the functionals $\wt J$ and
$\mathcal J$, one obtains
\be\label{3.5}
t^d \mathcal J(S_F) = \mathcal J(S_{F_{0,b_1t^2}}) + O(t^{d+1}).
\end{equation}

Our goal now is to evaluate the values of 
$\mathcal J(S_{F_{0,b_1t^2}})$. 
We first prove the following lemma estimating the distortion 
on geodesic disks in generic complex Banach manifolds $X$. 
Denote the \Ko \ and \Ca \ metrics on $X$ by $d_X$ and $c_X$ 
and assume that there is a \hol \ map $h: \D \to X$ such that   
$$
d_\D(t_1, t_2) = c_X(h(t_1), h(t_2)) = d_X(h(t_1), h(t_2)) 
$$ 
for any two distinct points $t_1, t_2 \in \D$.  
Then $h(\D)$ is a holomorphically embedded disk, geodesic
for both metrics $c_X$ and $d_X$. 
Such isometries between the hyperbolic metrics on $\D$ and $X$ 
are regarded as the {\em complex geodesics} (cf. \cite{Ve}).

\begin{lem} Let $X$ be a complex Banach manifold with
complete \ \Ko \ and \Ca \ metrics,
and let $h: \D \to X$ be a complex geodesic.
If the restriction of a  \hol \ map $X \to \D$ to the disk
$D = h(\D)$ has at $x_0 = h(0)$ zero of order $p \ge 1$,
then for all $t \in \D$,
\be\label{3.6}
|j \circ h(t)| \le \tanh d_X( h(t^p), \mathbf 0).
\end{equation}
The equality (even for one $t \ne 0$) only occurs for
defining functions of the \Ca \ distance
between the points $x_0$ and $x \in D$, then
$$
d_\D(j \circ h(t^p)), j \circ h(0)) = c_X(h_0(t^p), h(0))
\quad \text{for all} \ \ t \in \D.
$$
\end{lem}

\medskip
\noindent
\textbf{Proof}. First recall the Schwarz lemma for \sh \ functions
(we only need its simplest version).

\begin{lem} Let a function $u(t) : \ \D \to [0, 1)$ be
logarithmically \sh \ in the disk $\D$ and such that the ratio
$u(t)/|t|^m$ is bounded in a neighborhood of the origin for
some $m \ge 1$. Then
\be\label{3.7}
u(t) \le |t|^m \quad \text{for all} \ \ t \in \D
\end{equation}
and
\be\label{3.8}
\limsup\limits_{|t| \to 0} \frac{u(t)}{|t|^m} \le 1.
\end{equation}
Equality in (3.7), even for one $t_0 \ne 0$, or in (3.8), can
only hold for the function $u(t) = |t|^m$.
\end{lem}

The assumption on metrics $d_X$ and $c_X$ yields that the pluricomplex
Green function $g_X(x,y)$ of $X$ with a pole at $y$ cannot be equal 
identically $- \iy$ and that for any pair of points
$x, y \in X$ with equal \Ko \ and \Ca \ distances,
\be\label{3.9}
g_X(x, y) = \log \tanh d_X(x, y) = \log \tanh c_X(x, y);
\end{equation}
in addition,
$\lim g_X(x, y) = 0$ as $x$ tends to infinity (boundary of $X$),
for any fixed $y \in X$.
From (3.9) and Lemma 3.2,
$$
g_X(h(t^p), \mathbf 0) = \log |t|^p
\quad \text{for all} \ \ |t| <1.
$$

The function $j \circ h(t)|$ is logarithmically \sh \
on $\D$ and has zero of order $p$ at the origin,
hence by the same Schwarz's lemma,
$$
\log |j \circ h(t)| \le \log |t|^p,
$$
which implies the estimate (3.6), completing the proof of 
Lemma 3.1. 

\bk
We proceed to the proof of Theorem 1.1 and denote by $s$ 
the canonical complex parameter on the \Te \ disks 
in in $\mathcal F(\T)$ generated by admissible (that is,
nonvanishing on $\D^*$) functions
$$
F_{b_0,s}(z) = z + b_0 + s z^{-1}, 
$$ 
whose extremal extensions onto $\ov \D$ are the affine maps 
$z \mapsto z + b_0 + s \ov z$. 
These disks cover
$\D(S_{F_0,s}) \subset \T$.
If such $F_{b_0,s}$ is admissible only for $|s| < s_0 < 1$,
one can reparametrize it using the parameter $\sigma = s/s_0$ which
runs over the unit disk. For each $b_0$, 
\be\label{3.10}
d_{\T_1}(\mathbf 0, (S_{F_{b_0,s}}, b_0))
= d_\T(\mathbf 0, S_{F_{0,s}});
\end{equation}
in addition, every map
$$
\sigma \mapsto (S_{F_{b_0,\sigma}}, b_0 \sigma), \quad \sigma \in \T, 
$$
determines a complex geodesics in the space $\mathcal F(\T)$. 

Since by (2.14) the parameters $s$ and $t$ are
related near the origin by
$$
s = b_1 t^2  + O(t^3) \quad \text{as} \ \ t \to 0 
$$ 
and $b_1 \ne 0$, it follows from (3.5) that the restrictions 
of the functions (3.4) to any disk $\D(S_F{_{b_0,s}})$ 
have zero of order $d/2$ at the origin.  
Lemma 3.1 implies the bound  
$$
|g_m(S_{F_{b_0,s}})|
\le [\tanh d_\T(S_{F_{b_0,s}}, \mathbf 0)]^{d/2} = |s|^{d/2}, 
$$
or equivalently, 
\be\label{3.11}
|g_m(S_{F_{0,b_1t^2}})| \le |b_1|^{d/2} |t|^d + O(|t|^{d+1}),
\quad t \to 0.
\end{equation}
Note that by (3.3) the remainder in (3.11) does not depend on $m$.
Therefore,
\be\label{3.12}
\mathcal J(S_{F_{0,b_1t^2}}) = \sup_m |g_m(S_{F_{0,b_1t^2}})|
\le |b_1|^{d/2} |t|^d + O(|t|^{d+1}). 
\end{equation}
The relations (3.5) and (3.12) imply the equality 
\be\label{3.13}
\mathcal J(S_F) \le |b_1|^{d/2} 
= \left( \fc{|S_f(0)|}{6}\right)^{d/2} \le 1.
\end{equation}

Now observe that if the original functional $J(f)$ has an extremal
$f_0(z) = z + \sum\limits_2^\iy a_n^0 z^n$, for which the assumption
(1.10) is fulfilled, then the value of the left-hand side
in (3.21) on this function must be equal to $1$, because
in this case,
$$
\fc{|\wt J(S_{F_0}, a_2^0)|}{M(J)} = \mathcal J(S_{F_0}) = 1,
$$
thus
$$
\fc{1}{6} |S_{f_0}(0)| = |(a_2^0)^2 - a_3^0| = 1.
$$
As was mentioned, such equalities can only occur when $f_0$ either
is the Koebe function $\kp_\theta$ or it coincides with the odd
function $\kp_{2,\theta}$ defined by (1.6).
In addition, the extremality of $f_0$ implies
$$
|J(f_0)| = M(J) = \max \{|J(\kp_\theta)|, |J(\kp_{2,\theta})|\}.
$$

$(b)$ \ The functions $f \in S$ with $S_f(0) = 0$ omitted above
can be approximated (in $\B$-norm) by $f$ with $S_f(0) \ne 0$
by applying special \qc \ deformations of the plane given by the
following lemma from [Kr1, Ch. 4]. This lemma softens the
strongest rigidity of conformal maps.

\begin{lem} In a finitely connected domain $D \subset \hC$, let
there be selected a set $E$ of positive two-dimensional Lebesgue
measure and the distinct finite points $z_1, \dots, z_n$ with
assigned nonnegative integers $\a_1, \dots, \a_n$, respectively,
so that $\a_j = 0$ for $z_j \in E$.
Then, for sufficiently small $\ve > 0$ and $\ve \in (0, \ve_0)$,
for any given system of numbers
$\{w_{sj}\}, \ s = 0, 1, \dots, \a_j, \ \ j = 1, \dots, n$,
such that $w_{0j} \in D$,
$$
|w_{0j} - z_j| \le \ve, \ \ |w_{1j} - 1| \le \ve, \ \ |w_{sj}| \le \ve
\ \ (s = 2, \dots, \a_j, \ j = 1, \dots, n),
$$
there exists a \qc \ automorphism $h_\ve$ of the domain $D$, which is
conformal on the set $D \setminus E$ and satisfies
$h_\ve^{(s)}(z_j) = w_{sj}$ for all $s = 0, 1, \dots, \a_j$ and
$j = 1, \dots, n$, with dilatation $\|\mu_{h_\ve}\|_\iy \le M\ve$.
The constants $\ve_0$ and $M$ depend only on $D, \ E$ and the vectors
$(z_1, \dots, z_n), \ (\a_1, \dots, \a_n)$.

If the boundary $\G$ of domain $D$ is Jordan or belongs to the class
$C^{l,\a}$, where $0 < \a < 1$ and $l \ge 1$, one can take $z_j \in \G$
with $\a_j = 0$ or $\a_j \le l$, respectively.
\end{lem}

Now, let $f \in S^0$ have coefficients $a_2$ and $a_3$
related by $a_3 = a_2^2$, i.e., $b_1(f) := b_1(F_f) = 0$.
Since $f(\D^*)$ is a domain, one can take there
a set $E$ of positive measure and construct by
Lemma 3.5 for a sequence $\ve_n \to 0$ such variations
$h_n = h_{\ve_n}$ of $f$ that for each $n$,
$$
b_1(h_n \circ f)= b_1(f) + O(\ve_n) \ne 0, \quad
|J(h_n \circ f)| = |J(f)| + O(\ve_n) > |J(f)|.
$$
Since, by the previous step,
$$
|J(h_n \circ f)| \le \max \{|J(\kp_\theta)|, |J(\kp_{2,\theta})|\},
$$
the same estimate will hold also for $f$.

\bk
$(c)$ \ Finally, consider the case when $J$ has no extremals $f_0$
satisfying (1.10), and hence any extremal inversion $F_{f_0}$ is
of the form 
\be\label{3.14} 
F(z) = z + b_0 + b_m z^{-m} + b_{m+1} z^{-(m+1)} + \dots \
\quad (b_m \ne 0; \ |z| > 1)
\end{equation} 
with $m > 1$. 
If $m+1$ does not divide $d = d(J)$, we consider  the functional
$J^{m+1}$, which is $d(m + 1)$-homogeneous; otherwise one can use $J$.
 
Let us show that one can apply to $J^{m+1}$ the above arguments,
replacing $F_{0,b_1t^2}$ by the corresponding function 
$F_{m,t}$ represented by (2.9). Its Schwarzian relates to 
$S_{F_t}$ by  
$$
S_{F_t} = S_{F_{m,t}} + O(t^{m+1}), \quad t \to 0.
$$ 

If $m$ is odd, one immediately derives from Proposition 2.2 
(applied to $\mu_{F_{m,t}}(z) = t|z|^{m-1}/z^{m-1}$) that 
\be\label{3.15}
c_\T(\mathbf 0, S_{F_{m,t}}) = d_\T(\mathbf 0, S_{F_{m,t}})
= \tau_\T(\mathbf 0, S_{F_{m,t}})
= \tanh^{-1} \left( \fc{m + 1}{2} |b_m| |t|^{m+1} \right).   
\end{equation}  

If $m$ is even, we consider the map
$$
F_2(z) = F(z^2)^{1/2} = z + \fc{b_0}{2} \ \fc{1}{ z} +
\fc{b_m}{2} \ \fc{1}{z^{2m-1}} + \dots \ , 
$$
which is well defined, since $F(0) = 0$, and represents an odd 
function symmetric with respect to the origin.
Denote the Taylor coefficients of $F_2$ by $b_j^{(2)}$ and let
$\a_{mn}^{(2)}(F) = \a_{mn}(F_2)$.
Squaring 
$$ 
\mathcal R_2: \  F(z) \mapsto F(z^2)^{1/2}
$$
transforms the quadratic differentials
$\psi = \psi(z) dz^2$ in $\D$ into 
$\mathcal R_2^* \psi* = \psi_r(z^2) 4 z^2 dz^2$, 
which have zero of even order at the origin. 
Thus one can apply Proposition 2.2 to $F_2$, using instead of (2.8) 
the functions  
\be\label{3.16}
h_{2,\mathbf x}(\vp) =
\sum\limits_{p,q=1}^\iy \ \sqrt{p q} \ \a_{pq}^{(2)}(\vp) x_p x_p: \ 
\T \to \D. 
\end{equation} 

Indeed, all maps $f$ and $F_f$ are completely normalized
and their \Be \ coefficients and Schwarzians are related by
$$
\mu_F(z) = \mu_f(1/z) z^2/\ov z^2, \quad
S_F(z) = - S_f(1/z) z^{-2}. 
$$
Applying the Cauchy formula for derivatives of
\hol \ functions, one derives that the Taylor coefficients 
$a_n, \ n \ge 2$,  
and $b_j, \ j \ge 0$, depend holomorphically on 
\Be \ coefficients $\mu_F \in \Belt(\D)_1$ and on 
Schwarzians $S_F \in \T$.

On the other hand, each of the coefficients $b_j^{(2)}$ and 
$\a_{pq}^{(2)}$
of $F_2$ is represented as a polynomial of a finite number of
initial coefficients $b_0, b_1, \dots \ , b_s$
of the original function $F$ (noting that the free term $b_0$ is  
uniquely determined by assumption $F(0) = 0$).
Thus the maps (3.16) also depend holomorphically on $\mu_F$
and $S_F$. 
Finally, the transform $\mathcal R_2$ preserves 
\qc \ dilatations.  

These properties imply that $S_{F_{m,t}}$ satisfies 
all equalities in (3.15) and hence ranges in a geodesic 
disk in $\T$.  

Combining with (2.14), one now obtains instead of (3.13) the bound
$$
\mathcal J(S_F)^{m+1}
\le \left(\fc{m + 1}{2} |b_m|\right)^d \le 1, 
$$
or
\be\label{3.17}
\mathcal J(S_F) \le \left(\fc{m + 1}{2} |b_m|\right)^{d/(m+1)}
\le 1.
\end{equation}

To examine the case of equality in (3.17), observe that since
$\mathcal J(S_{F_{f_0}}) = 1$ for any extremal $f_0$, (3.17) yields
(denoting the coefficients of $F_{f_0}$ by $b_j^0$)
$$
\fc{m + 1}{2} |b_m^0| = 1.
$$
The functions (3.13) with $m > 1$ satisfy
$b_1 = \dots = b_{m-1} = 0$, which allows us to estimate $b_m$
using the well-known inequalities of Golusin and Jenkins (see
[Go, Ch. XI], \cite{Je}). These inequalities hold for a more general
univalent $F(z) = z + \sum\limits_0^\iy b_n z^{-n}, \ |z| > 1$
(but not for all $F \in \Sigma$), and in our special case imply
the bound
$$
|b_m| \le 2/(m + 1).
$$
Moreover, the equality here only occurs for the function
(2.9) with $|t| = 1$
(and its admissible translations $F_{m,t}(z) + c$), which
also implies
$M(J) = |J(\kp_{m,\theta})|$.

We have established that any extremal function $f_0$ maximizing
$|J(f)|$ must be of the form (1.8). The theorem is proved.

\bk 
\noindent 
\textbf{Remark}. Lemma 3.1 is interesting for its own sake and has 
other applications. Thus we also present a different proof which
does not involve Green's function and Lemma 3.2.
By assumption,
$$
g \circ h(t) = c_p t^p + c_{n+1} t^{p+1} + \dots \ (c_p \ne 0).
$$
Let $\{g_m\}$ be a maximizing sequence for the \Ca \ distance
$c_X(x_0, x)$ with $g_m(x_0) = 0$, where $x \in D = h(\D)$ is
distinct from $x_0$, i.e.,
$$
c_X(x_0, x) = \lim\limits_{m\to \iy} d_\D(0, g_m(x)).
$$
The restrictions of $g_m$ to $D$ are convergent locally
uniformly on $D$ to a \hol \ function $g_0$ on this disk
with $g_0(x_0) = 0$ so that
$d_\D(0, g_0(x)) = c_X(x_0, x), \ x \in D$. 
Since both metrics $c_X$ and $ d_X$ are hyperbolic isometries
between $\D$ and $D$, the map
$g_0 \circ h$ is the identity on $\D$; hence,
$g_0^p \circ h(t) = t^p$ for all $t \in \D$, and
$$
|j \circ h(t)| \le |g_0^p \circ h(t)|. 
\quad t \in \D.
$$
Together with the equality 
$g_0 \circ h(t^p) = t^p, \ t \in \D$, this implies (3.6).
The case of equality in (3.6) follows from Schwarz's
lemma for \hol \ functions with zero of a prescribed
order at the origin.

\section{Proof of Theorem 1.2}

Note that from (1.3),
$$
a_n^2 - a_{2n-1} = b_1 b_0^{2n-4} +
\text{lower terms with respect to} \ b_0.
$$
We have to show that for all $m > 1$,
\be\label{4.1}
|J_n(\kp_{m,\theta})| < J_n(\kp),
\end{equation}
then Theorem 1.1 implies that only the Koebe function
is extremal for Zalcman's functional.

This inequality is trivial for $m = 2$, because
the series (1.7) yields
$$
|J_n(\kp_{2,\theta})| \le 2 < J_n(\kp).
$$
For $m \ge 3$, we apply a result of \cite{Kr5} solving
the coefficient problem for univalent functions
with \qc \ extensions with small dilatations.
Denote by $S_\iy(k)$ the subclass of $S^0$ consisting of
the functions $f \in S$ having $k^\prime$-\qc \ extensions
$\wh f$ to $\hC \ (k^\prime \le k)$, which satisfy
$\wh f(\iy) = \iy$, and let
$$
f_{1,t}(z) = \fc{z}{(1 - kt z)^2}, \quad |z| < 1, \ \ |t| = 1.
$$

\begin{prop} \cite{Kr5} For all
$f(z) = z + \sum\limits_2^\iy a_n z^n \in S_\iy(k)$
and all $k \le 1/(n^2 + 1)$, we have the sharp bound
\be\label{4.2}
|a_n| \le \fc{2 k}{n - 1},
\end{equation}
with equality only for the functions
\be\label{4.3}
f_{n-1,t}(z) = f_{1,t}(z^{n-1})^{1/(n-1)}
= z + \fc{2 k t}{n - 1} z^n + \dots, \quad n = 3, 4, \dots
\end{equation}
\end{prop}

Note that every function (4.3) admits a \qc \ extension onto
$\D^*$ with \Be \ coefficient
$\mu_n(z) = t |z|^{n+1}/z^{n+1}$ 
and $\wh f_{n-1,t}(\iy) = \iy$. Accordingly,
$F_{n-1,t}(z) = 1/f_{n-1,t}(1/z) \in \Sigma^0$ admits a \qc \
extension onto the unit disk with $F_{n-1,t}(0) = 0$ and
$\mu_{F_{n-1,t}}(z) = t |z|^{n-1}/z^{n-1}$ for $|z| < 1$. 
Another essential point is that for any function
$$
F_{n-1}(z) := F_{n-1,1}(z) = 1/\kp(1/z^{n-1})^{1/(n-1)},
$$
its homotopy disk $\D(S_{F_{n-1}})$ in $\T$ is of \Te \ type.
Together with estimate (4.2), this implies that for any $m > 2$
and small $r > 0$,
$$
|J_n(\kp_{m,r})| < r (n - 1 )^2;
$$
thus
$$
|J_n(\kp_{m,\theta})| < (n - 1 )^2 = J_n(\kp),
$$
completing the proof of Theorem 1.2.

\bk
\noindent
\textbf{Remark}. The above arguments work well
also in the case of functionals obtained by suitable
perturbation of $J_n(f)$. For example, one can take
$$
J(f) = a_n^2 - a_{2n-1} + P(a_3, \dots, \ a_{2n-2}),
$$
where $P$ is a homogeneous polynomial of degree $2n-2$,
$$
P(a_3, \dots, \ a_{2n-2}) = \sum\limits_{|k|=2n-2}
c_{k_3, \dots, k_n} a_3^{k_2} \dots a_n^{k_{2n-2}},
$$
and $|k| := k_3 + \dots \ + k_{2n-2}, \ a_j = a_j(f)$,
assuming that this polynomial has nonnegative coefficients
and satisfies
$$
\max_S |P(a_3, \dots, \ a_{2n-2})| < \fc{(n - 1)^2}{2}.
$$
For any such functional, only the Koebe function
is extremal.

\section{Some new distortion theorems for higher coefficients}

As was mentioned, Theorem 1.1 provides various new
distortion estimates. For example, we obtain the following
generalization of the inequality
$|a_2^2 - a_3| \le 1$ to higher coefficients.

\begin{thm} For all $f \in S$ and integers $n > 3$
and $p \ge 1$,
$$
|a_n^p - a_2^{p(n-1)}| \le 2^{p(n-1)} - n^p.
$$
This bound is sharp, and the equality only occurs for the Koebe
function $\kp_\theta$.
\end{thm}

\medskip
\noindent
\textbf{Proof}. Since $b_0 = - a_2$, the relation (1.3) yields
$$
I_n(f) := a_n - a_2^{n-1}
= (n - 2)(- 1)^{n-1} b_1 b_0^{n-3} +
\text{lower terms with respect to} \ \ b_0.
$$
This functional satisfies the assumptions of Theorem 1.1.
The same arguments as in the proof of Theorem 1.2 imply
$$
|I_n(\kp_{m,\theta})| < |I_n(\kp_\theta)| \quad \text{for all} \ \
m \ge 2,
$$
completing the proof.

In the same way, one obtains

\begin{thm} For all $f \in S$ and integers $n > 2$ and $p \ge 1$,
$$
|a_{n+1}^p - a_2^p a_n^p| \le 2^p n^p - (n + 1)^p,
$$
with equality only for $f = \kp_\theta$.
\end{thm}

\medskip
{\small\em{
\leftline{Department of Mathematics, Bar-Ilan University}
\leftline{52900 Ramat-Gan, Israel}
\leftline{and Department of Mathematics, University of Virginia,}
\leftline{Charlottesville, VA 22904-4137, USA}}}

\end{document}